\documentclass[11pt,twoside]{article}
\usepackage{graphicx}
\usepackage{amssymb}
\usepackage{epsfig}
\usepackage{amsmath}
\usepackage{amsfonts}
\usepackage{color}
\newtheorem{theorem}{\bf Theorem}[section]

\newtheorem{remark}{\bf Remark}[section]
\newtheorem{lemma}{\bf Lemma}[section]
\newtheorem{definition}{\bf Definition}[section]

\newsavebox{\savepar}

\begin{document}
\title{\bf Approximate Controllability of a Class of Partial  Integro--Differential Equations of Parabolic Type}
\author{{\scshape Anil Kumar}\footnote{Email:
anil@goa.bits-pilani.ac.in}, {\scshape Amiya K. Pani}\footnote{Email:
akp@math.iitb.ac.in} \hspace{.05cm} {\scshape and} {\scshape Mohan C. Joshi}\footnote{Email: mcj@math.iitb.ac.in}\\
$^{*}$ Department of Mathematics, BITS Pilani KK Birla Goa Campus \\
Zuarinagar, Goa-403726, India.\\
$~^{\dag}$  Department of Mathematics, Industrial Mathematics Group\\  Indian Institute of Technology Bombay,
Powai, Mumbai-400076, India.\\
$~{^\ddag}$ Indian Institute of Technology Gandhinagar\\
Ahmedabad, Gujarat-382424, India.}
\maketitle
\pagenumbering{arabic}

\begin{abstract}
In this paper, we discuss the distributed control problem governed by the following parabolic
integro-differential equation (PIDE) in the abstract form
\begin{eqnarray*}
\frac{\partial y}{\partial t} + A y &=& \int_0^t B(t, s) y(s) ds +
Gu, \;\; t \in [0, T], \;\;\;\;\;\;\;\;\;\;\;\;\;\;\; \hfill{(\ast)}\\
y(0) &=& y_0 \; \in X, \nonumber
\end{eqnarray*}
where, $y$ denotes the state space variable, $u$ is the control variable, $A$ is a self adjoint, positive definite linear (not necessarily bounded) operator in a Hilbert space $X$ with dense domain $D(A) \subset X,$ $B(t,s)$ is an unbounded operator, smooth with respect to $t$ and $s$ with $D(A) \subset D(B(t,s))
\subset X$ for $0 \leq s \leq t \leq T$ and $G$ is a bounded linear operator from the control space to $X.$  Assuming that the corresponding  evolution equation ($B \equiv 0$ in ($\ast$))
is approximately controllable, it is shown that the set of approximate controls of 
the distributed control problem ($\ast$) is nonempty.
The problem is first viewed as constrained optimal control problem and then it is approximated by  unconstrained problem with a suitable penalty function. The optimal pair of the constrained problem is obtained as the limit of optimal pair sequence  of the unconstrained problem. The approximation theorems, which guarantee the convergence of the numerical scheme to the optimal pair sequence, are also proved.\\

\end{abstract}

\section{Introduction}
Consider the following  parabolic integro-differential equation with distributed control
\begin{eqnarray}\label{e1.1}
\frac{\partial y}{\partial t} + A y(t) &=& \int_0^t B(t, s) y(s) ds +
Gu(t), \;\; t \in [0, T], \\
y(0) &=& y_0 \; \in X, \nonumber
\end{eqnarray}
where $X$ denotes a real Hilbert space, $y$ is a state variable, $u$ represents a control variable,
$A$ is a self adjoint, positive definite linear operator in $X$ with dense domain $D(A) \subset X,$
$B(t,s)$ is also a linear and unbounded operator with $D(A) \subset D(B(t,s)) \subset X$ for $0 \leq s
\leq t \leq T$ and $G$ is a bounded linear map from the control space $U$ to $X.$

For an example, let $\Omega$ be a bounded domain in $\mathbb{R}^d$ with
smooth boundary $\partial \Omega$. For fixed $T < \infty$, let $Q
= (0, T) \times \Omega$ and $\Sigma = (0, T) \times \partial
\Omega$. Further, set $A$ as a second-order linear self-adjoint elliptic partial differential operator defined by
\begin{eqnarray}\label{ee1.2}
A = - \sum_{i, j = 1}^d \frac{\partial}{\partial x_j}
\left(a_{i j}(x) \frac{\partial}{\partial x_i} \right) + a_0(x) I,
\end{eqnarray}
where the matrix $(a_{i j}(x))$ is symmetric and positive definite,
$a_0 \geq 0$ on $\bar{\Omega}$ and  $B(t, s)$ is a general
second-order partial differential operator of the form
\begin{eqnarray}\label{ee1.3}
B(t, s) = - \sum_{i,j=1}^d \frac{\partial}{\partial x_j}
\left( b_{i,j}(t, s; x) \frac{\partial}{\partial x_i} \right) +
\sum_{j=1}^d b_j(t,s;x) \frac{\partial}{\partial x_j} + b_0(t,s;x)
I,
\end{eqnarray}
with smooth coefficients $b_{i,j}, b_j$ and $b_0$. Let $X = L^2(\Omega), \; D(A) =
H^2(\Omega) \cap H_0^1(\Omega)$
and $D(B) = H^2(\Omega)$, then the following problem
\begin{eqnarray}\label{e1.2}
\frac{\partial y(t,x)}{\partial t} + A(x) y(t,x) &=& \int_0^t B(t,
s) y(s) \;ds + G u(t,x) \; \mbox{in} \; Q, \nonumber\\
y(t,x) &=& 0 \; \mbox{on} \; \Sigma, \\
y(0,x) &=& y_0(x) \; \mbox{in}\; \Omega, \nonumber
\end{eqnarray}
becomes a particular case of the abstract problem (\ref{e1.1}), where $y_0 \in X$.

Parabolic integro-differential equations of the type (\ref{e1.1}) occur in many applications such as 
heat conduction in materials with memory, compression of poroviscoelastic media, nuclear reactor dynamics,
etc. (see, Cushman \emph{et al.} \cite{C1}, Dagan \cite{D1}, Renardy \emph{et al.} \cite{R1}).

For control  problems of the heat equation  with memory, that is, when $A= -\Delta $ and 
$B(t,s)= - a(t-s) \Delta u (s)$ in (\ref{e1.1}), where $a(\cdot)$ is a completely monotone convolution kernel, 
Barbu and Iannelli \cite{BI-2000} have discussed approximate controllability using Carleman estimates. Later on,
Pandolfi \cite{P-2009} has considered  Dirichlet boundary controllability of heat equation with memory in one 
space dimension by employing Riesz systems. There are several negative results like lack of controllability
of such systems, see, \cite{IP-2009}, \cite{HP-2012} and \cite{GI-2013}.  Subsequently, Fu {\it et al.} 
\cite{FYZ-2009} have established controllability and observability  results for 
a heat equation with hyperbolic memory  kernel under  general  geometric conditions and using 
Carleman estimates. Further, Pandolfi \cite{P-2005} has employed cosine operator approach to discuss the
exact controllability  results for the Dirichlet boundary control of  the Gurtin-Pipkin model 
which displays a hyperbolic behaviour. On second order integro-differential equations, Kim \cite{K1, K2}
has established reachability results using continuation arguments and multiplier techniques combined with
compactness property. Wang  and Wei \cite{W1} have proved some sufficient conditions for the 
controllability of parabolic integro-differential systems in a 
Banach space. A result in the direction of approximate controllability of integro-differential equations (IDE) 
using Carleman estimates and continuation argument has been proved by Lefter {\em et. al.} in \cite{L1}. Loreti
and Sfrorza \cite{LS-2010} have analyzed  reachability problems  for a class of  IDE using Hilbert uniqueness
results. In the present article, an attempt has been made to discuss approximate controls of a distributed 
control problem for  a general class of  partial integro-differential equations of parabolic 
type (\ref{e1.1}), under the assumption that the corresponding
parabolic equation is approximately controllable. Firstly, the control problem is viewed as an optimal 
control problem, and using operator theoretic form, an optimal pair of solution is derived, which, in turn,
provides a proof for the approximate controllability.The present proof is constructive in its approach 
and avoids of using Carleman  estimates and continuation of argument, etc. Finally, some approximate theorems
are established and numerical experiments using finite element method are conducted to  confirm our theoretical
findings.

Numerical solution by means of finite element methods has been invested by several authors, when $u$ is a given function and $G=I.$  In \cite{T1}, Thom\'ee {\em et. al.} have considered backward Euler methods and obtained related error estimates for non-smooth data. Pani {\em et. al.} \cite{P1} have used energy arguments and the duality technique to obtain error estimates for time dependent parabolic integro-differential equations with smooth and non-smooth initial conditions. 
Lasiecka \cite{LL1, LL2} have considered  optimal control problems for 
linear parabolic equations, which are approximated by a semidiscrete finite element method or 
Ritz-Galerkin scheme and then the convergence of optimal controls are derived. Moreover,
 Shen {\em et. al.} \cite{S1} have developed the finite element and backward Euler scheme for the space
and time approximation of a constrained optimal control problem governed by a parabolic integro-differential 
equation. Further, in \cite{W2} Shen {\em et. al.} have  discussed mathematical formulation 
and optimality conditions for a quadratic optimal control problems for a quasi-linear integral differential 
equation and some {\it a prior} error estimates are also discussed.

In order to motivate our main results, we first define the operator $\tilde{B}$ as
$$ (\tilde{B} y)(t) = \int_0^t B(t, \tau) y(\tau) d\tau. $$

Since $A$ generates a $C_0$-semigroup $\{S(t)\}_{t \geq 0}$ of bounded linear operators on $X,$ then for a given $u \in U$
and $y_0 \in D(A)$, the mild solution for the system (\ref{e1.1}) is given by
\begin{equation}\label{e1.3}
y(t) = S(t) y_0 + \int_0^t S(t - \tau) \tilde{B} y (\tau) \;d\tau + \int_0^t S(t - \tau) Gu(\tau) \;d\tau
\end{equation}
(refer, Pazy \cite{P}). This correspondence which assigns a unique $y \in Z = L^2(0, T; X)$ to a given $u \in U$, will be denoted by a solution operator, say $W$ {\em i.e.} $Wu = y$. Also, set $Y= L^2(0,T;U)$.

The system (\ref{e1.1}) is said to be approximately controllable if for given functions $y_0,$ $\hat{y} \in X$ and a $\delta>0$, there exists a control $u \in U$ such that the corresponding solution $y$ of system (\ref{e1.1}) also satisfies $\|y(T)-\hat{y}\|_{X}\leq \delta.$

In view of (\ref{e1.3}), for such control $u$, we arrive at
\begin{equation}\label{e1.4}
\hat{y} = S(T) y_0 + \int_0^T S(T - \tau) \tilde{B} y (\tau)
\;d\tau + \int_0^T S(T - \tau) Gu(\tau) d\tau,
\end{equation}
where $\hat{y} = y(T)$.
Setting the operator $L : U \rightarrow X$ as
\begin{equation}\label{15}
L u = \int_0^T S(T - \tau)u(\tau) d\tau,
\end{equation}
then the last term on the right hand side of  (\ref{e1.4}) becomes $LGu.$
Now, the adjoint operator $L^* : X  \rightarrow Z$  of $L$ becomes
$$(L^*z)(\tau) =S(T-\tau)z, \; \tau \in [0,T] \;\mbox{and} \; z \in X.$$

If $G^*$ is the adjoint operator of the operator $G$, then it follows that
$$( G^* L^*z)(\tau) =  G^* S(T-\tau)z, \; \tau \in [0,T] \;\mbox{and} \; z \in X.$$

Thus, the equation (\ref{e1.4}) can be written equivalently as an operator equation
\begin{equation}\label{e1.5}
\hat{z} = L \tilde{B} y + LGu,
\end{equation}
where $\hat{z} = \hat{y} - S(T) y_0$.

Define for $\delta >0,$ the set $U_{\delta} \subset Y$ of admissible controls of (\ref{e1.1}) by
$$ U_{\delta} = \left\{ u \in Y: \; \|L \tilde{B} y + LGu - \hat{z} \|_{X} \leq \delta \right\}. $$

It is a closed, convex and bounded (possibly empty) subset of $Y.$

\begin{definition}
The problem (\ref{e1.1}) is approximately controllable if for every $y_0, \, \hat{y} \in X$ and $\delta>0,$
there exists $u \in Y$ such that $U_{\delta} \neq \emptyset.$
\end{definition}

We now define our main problem as

\noindent
{\bf Main Problem.} Find
\begin{itemize}
 \item [(i)] if $U_{\delta} \neq \emptyset$ for each $\delta > 0 $ and 
 \item [(ii)] if so  determine  $u_{\delta}^* \in U_{\delta}$ such that
\begin{equation}\label{e1.6}
    J(u_{\delta}^*) = \inf_{u \in U_{\delta}} J(u)
\end{equation}
where $J(u) =  \frac{1}{2}\|u\|_Y^2.$
\end{itemize}

\begin{definition}
For a given $\delta > 0$, let  $u_{\delta}^* \in U_{\delta}$ is a solution of the problem (\ref{e1.6}) with
$y_{\delta}^* \in X$ as the corresponding mild solution of the system (\ref{e1.1}),  then the pair
$(u_{\delta}^*, y_{\delta}^*)$ is called optimal pair of the constrained optimal control problem (\ref{e1.6}).
\end{definition}

Our main thrust is to establish the existence of the optimal pair
$(u_{\delta}^*, y_{\delta}^*)$ of the constrained optimal control problem (\ref{e1.6}) and thereafter  present a numerical scheme for approximating the optimal pair.  Under the assumption $B \equiv 0$, in section 2 we first show that the set $U_{\delta}$ of admissible controls is nonempty. Then the optimal pair
$(u_{\delta}^*, y_{\delta}^*)$ is obtained as a limit of the sequence of an optimal pair $(u_{\epsilon}^*, y_{\epsilon}^*)$, where $u_{\epsilon}^*$ minimizes the unconstrained functional
$J_{\epsilon}(u)$ over the whole space $Y$ defined by
\begin{equation}\label{e1.7}
J_{\epsilon}(u) = J(u) + \frac{1}{2\epsilon} \left\|L u + L \tilde{B} W u - \hat{z} \right\|_X^2,
\end{equation}
where $W$ is the operator which assigns to each control $u_{\epsilon}^*$ the solution  $y_{\epsilon}^*$ of (\ref{e1.1}). We shall refer to $(u_{\epsilon}^*, y_{\epsilon}^*)$ as the
optimal pair corresponding to the unconstrained problem.


The plan of this paper is as follows: In Section 2, we have shown that the set of admissible control $U_{\delta}$ is nonempty under the assumption that the corresponding linear system is approximately controllable. The optimal pair $(u_{\delta}^*, y_{\delta}^*)$ of the constrained problem (\ref{e1.6}) is obtained as a limit of the optimal pair sequence $(u_{\epsilon}^*, y_{\epsilon}^*)$, where $u_{\epsilon}^*$ minimizes the unconstrained functional $J_{\epsilon}(u)$ defined by (\ref{e1.7}). We present approximation theorems which guarantee the convergence of the numerical scheme to the optimal pair in Section 3. Error estimates are derived for the final state of the problem in section 4 with an application. In section 5, we conclude this paper by providing some numerical experiments to demonstrate the applicability of our results.

\noindent
\section{Existence of optimal control and convergence to the control problem}
\setcounter{equation}{0}
In this section, we first show that the set $U_{\delta}$ of admissible controls is nonempty. Here, we first
make the following assumptions for the problem (\ref{e1.1}):
\begin{description}
\item{(A1)}
The set $\{S(t)\}_{t \geq 0}$ of $C_0$-semigroup of bounded linear operators on $X$, generated by $(-A)$ is uniformly bounded, that is, there exists $\beta > 0$ such that $ \| S(t)\|_X \leq \beta$, for all $t \in [0,T]$.
\item{(A2)}
The operator $B(t,\tau)$ is dominated by $A$ together with certain derivatives with respect to $t$ and $\tau$, that is,
$\|A^{-1} B(t,\tau) \phi \| \leq \alpha \|\phi\| \;\;  \forall \;\; \phi \in D(B(t,\tau)), \;\; 0 \leq
\tau \leq t \leq T.$ 
\item{(A3)}
The system (\ref{e1.1}) with $B \equiv 0$ is approximately controllable.
\item{(A4)}
The operator $G: L^2(0,t;U) \rightarrow L^2(0,T;X)$ is a bounded linear operator.
\end{description}

The following lemma is related to the assumption  $(A3).$
\begin{lemma}\label{l3.1}
The system (\ref{e1.1}) with $B \equiv 0$ is approximately controllable on $[0, T]$ if and only if one of the following statement holds:
\begin{description}
\item[(i)] $\overline{Range(LG)} = X$.
\item[(ii)] $Kernel(G^*L^*)= \{0\}$.
\item[(iii)] For all $z \in X$, there holds for $\delta \in (0,1)$
$$LGu_{\delta}= z-\delta\Big(\delta I + LG G^* L^*\Big)^{-1} z, $$
where
$u_{\delta}:=G^*L^*\Big(\delta I + LG \;G^* L^*\Big)^{-1} z.$
\item [(iv)] $\displaystyle \lim_{\delta \to 0^{+}} \delta \Big(\delta I + LG \;G^* L^*\Big)^{-1} z =0.$
\end{description}
\end{lemma}
For a proof, we refer to  Curtain {\em et. al.} \cite{c2, c3}.

As a consequence, it is observed that
$$\lim_{\delta \to 0^{+}} LG u_{\delta}= z$$
and the error $e_{\delta} z$ due to this approximation is given by
$$e_{\delta} z =\delta \Big(\delta I + LG G^* L^*\Big)^{-1} z =0.$$

For approximate controllability of the problem (\ref{e1.1}), we rewrite its  controllability equation as
\begin{equation}\label {eq:control}
u_{\delta} := G^*L^* \Big(\delta I + LG G^*L^*\Big)^{-1} (\hat{z}- L\tilde{B} y),
\end{equation}
where $\tilde{z}= y(T)- S(T) y_0.$

Now for a fixed $z \in Z$, consider the following linear parabolic integro--differential system which is
indexed by $z$
\begin{eqnarray}\label{e3.2}
\frac{\partial y_z}{\partial t} + A y_z &=& \int_0^t B(t, \tau) z(\tau) d\tau + Gu_z, \;\; t \in [0, T], \\
y_z(0) &=& y_0 \in X. \nonumber
\end{eqnarray}
The mild solution $y_z \in Z$ of the above system is given by
\begin{eqnarray}\label{e3.3}
y_z(t) = S(t) y_0 + \int_0^t S(t - \tau) \tilde{B} z(\tau) \,
d\tau + \int_0^t S(t - \tau) Gu_z(\tau) d\tau,
\end{eqnarray}
and hence,
\begin{eqnarray} \label{e3.4}
\tilde{y}_z \equiv  y_z(T) - S(T) y_0 = \int_0^T S(T - \tau) \tilde{B}
z(\tau) \, d\tau +  \int_0^T S(T - \tau) Gu_z(\tau) \, d\tau.
\end{eqnarray}

In operator theoretic form, (\ref{e3.4}) reduces to the operator equation
\begin{eqnarray} \label{e3.5}
LGu_z = \hat{y}_z - L \tilde{B} z,
\end{eqnarray}
for each fixed $z$.

Now under the assumption $(A3),$ the system (\ref{e1.1}) with $B \equiv 0$ is approximately controllable, and
hence, Lemma \ref{l3.1} implies that $(\delta I +G^* L^* L G)$ is boundedly invertible. For
approximate controllability of (\ref{e3.2}), we observe that for a given state
$\hat{z} = \hat{y} - S(T) y_0 \in X$, and for $\delta >0,$ a control $u_{\delta,z}$ solves
\begin{eqnarray} \label{e3.6}
u_{\delta,z} = G^*L^* (\delta I + LG\;G^* L^*)^{-1} \left[ \hat{z} - L \tilde{B} z \right].
\end{eqnarray}
To keep the notation simple and where there is no confusion, we write $u_{\delta,z}$ simply by $u_z.$

Denote the operator $\mathcal{M} z:=G^*L^* (\delta I + LG G^* L^*)^{-1} \left[ \hat{z} - L \tilde{B} z \right]$
and consider for $\delta \in (0,1]$ the family of operators $R_{\delta} : Z \rightarrow Z$,
which assigns a solution $y_z$ of (\ref{e1.1}) (given by
(\ref{e3.3})), corresponding to $z \in Z,$ that is,
\begin{eqnarray}\label{e3.7}
R_{\delta} z(t) &=& S(t) y_0 + \int_0^t S(t - \tau) \left(\tilde{B} z(\tau) + LG {\mathcal{M}} z(\tau) \right) d\tau.
\end{eqnarray}

Define the operator $K$ by
\begin{equation}\label{e3.8}
(K y)(t) = \int_0^t S(t - \tau) y(\tau) d\tau .
\end{equation}

Rewriting equation (\ref{e3.7}) in operator form as
\begin{equation}\label{e3.77}
R_{\delta} z(t) = S(t) y_0 + K \tilde{B}z(t) + KG \mathcal{M} z(t).
\end{equation}

First of all, we need to prove that for each fixed $\delta \in (0,1]$ the operator $R_{\delta}$ has
fixed point, say, $z_{\delta}.$

Now, the  following lemmas deals with some properties of $K,$ and $L.$
\begin{lemma}\label{l2}
Let the assumptions $(A1)$ and $(A2)$ be satisfied and let the operators $L : U \rightarrow X$ and $K : Z \rightarrow Z$ be defined by
(\ref{15}) and (\ref{e3.8}) respectively. Then the following estimates hold
$$ \| (K \tilde{B} y)(t) \|_X \leq C \int_0^t \|y(s)\|_X ds, $$
and
$$ \| L \tilde{B} y \|_X \leq C \int_0^T \|y(s)\|_X ds .$$
\end{lemma}
{\bf Proof.} From the definition of $K$ and $\tilde{B}$, we rewrite using semigroup property to get
\begin{eqnarray*}
(K \tilde{B} y)(t) &=& \int_0^t S(t-\tau) \int_0^{\tau} B(\tau,
s) y(s) ds d\tau \\
&=& \int_0^t S(t-\tau) A A^{-1} \int_0^{\tau} B(\tau,
s) y(s) ds d\tau \\
&=& - \int_0^t \frac{d}{d \tau} S(t-\tau) \int_0^{\tau}A^{-1}
B(\tau, s) y(s) ds d\tau.
\end{eqnarray*}
After integrations by parts, we arrive at
\begin{eqnarray*}
(K \tilde{B} y)(t) &=& \left[ S(t-\tau)\int_0^{\tau} A^{-1} B(\tau,
s) y(s) ds \right]_0^t - \int_0^t S(t-\tau) \int_0^{\tau} A^{-1}
B_{\tau}(\tau, s) y(s) ds d\tau \\
&&~~~~~~~~~~~~~~~~~~~~~ - \int_0^t S(t-\tau)A^{-1} B_{\tau}(\tau,
\tau) y(\tau) d\tau \\
&=& \int_0^t A^{-1}B(t, s) y(s) ds - \int_0^t S(t-\tau)
\int_0^{\tau} A^{-1} B_{\tau}(\tau, s) y(s) ds d\tau \\
&&~~~~~~~~~~~~~~~~~~~~~ - \int_0^t S(t-\tau)A^{-1} B_{\tau}(\tau,
\tau) y(\tau) d\tau.
\end{eqnarray*}
We note that
\begin{eqnarray*}
\|(K \tilde{B} y)(t)\|_X &\leq& \int_0^t \|A^{-1}B(t, s) y(s)\| ds +
\int_0^t \|S(t-\tau)\| \int_0^{\tau} \|A^{-1} B_{\tau}(\tau, s)
y(s)\| ds d\tau \\
&&~~~~~~~~~~~~~~~~~~~~~ + \int_0^t \|S(t-\tau)\| \, \|A^{-1}
B_{\tau}(\tau, \tau) y(\tau)\|d\tau,
\end{eqnarray*}
and hence,
\begin{eqnarray*}
\|(K \tilde{B} y)(t)\|_X &\leq& \alpha \int_0^t \|y(s)\|_X ds +
\alpha \|A^{-1}\| (1+\beta) \int_0^t  \|y(\tau)\|_X  d\tau
+ \alpha \beta \int_0^t \| y(\tau)\|_X d\tau\\
&\leq& \alpha (1 + \beta) (1+ \|A^{-1}\|) \int_0^t \|y(\tau)\|_X d\tau.
\end{eqnarray*}
Thus, we now arrive at
$$ \|(K \tilde{B} y)(t)\|_X \leq C \int_0^t \|y(\tau)\|_X d\tau, $$
where $C$ is a generic constant which depend on $\alpha, \beta$ and
$\|A^{-1}\|$. Similarly we can show by using the definition of $L$ and
$\tilde{B}$ that
$$ \| L \tilde{B} y \|_X \leq C \int_0^T \|y(\tau)\|_X d\tau .$$
This completes the proof of the lemma. \hfill{$\blacksquare$}

On the lines of Lemma \ref{l2}, we have the following result.
\begin{lemma}\label{l33}
Under the assumptions $(A1), (A2)$ and $(A4),$ the following estimate holds
$$ \left\|\Big(KLG\mathcal{M}y\Big)(t)\right\|_X \leq C_1 \left(\|\hat{z}\| + C \int_0^t \|y(\tau)\|_X d\tau\right).$$
where $C_1$ depends on $T, \beta, \|LG\|,  \|G^*L^*\|$ and $\|(\delta I + LG G^* L^*)^{-1}\|$.
\end{lemma}

A variation of Banach contraction mapping principle will help in the proof of the following theorem, which provides the approximate controllability of the system (\ref{e1.1}).
\begin{theorem}\label{t11}
Under the assumption $(A1)-(A4)$, the operator $R_{\delta}^n$ is a contraction on the space $Z$ for some positive integer $n$. Moreover, for any arbitrary
$z_0\in X$, the sequence of iterates $\{z_{\delta,k}\}$, defined by
\begin{equation}\label{e3.9}
    z_{\delta,k+1} = R_{\delta}^n z_{\delta,k}, \;\;\; k = 0, 1, 2, \ldots
\end{equation}
with $z_{\delta,0}=y_0$
converges to $y_{\delta}^*$, which is a mild solution of the system
(\ref{e1.1}). Further, $u_{\delta,k} = \mathcal{M} z_{\delta,k}$ is such that $u_{\delta,k}$ converges to
$u^*_{\delta} = \mathcal{M} y^*_{\delta},$ and the system  (\ref{e1.1})  is approximately controllable.
\end{theorem}
{\bf Proof.} Let $z_1, z_2 \in Z$, then use of (\ref{e3.77}) yields
\begin{eqnarray*}
  (R_{\delta} z_1 - R_{\delta} z_2)(t) &=& K\tilde{B}(z_1(t)-z_2(t))+KG \mathcal{M}(z_1(t)-z_2(t)).
\end{eqnarray*}
Using Lemma \ref{l2} and \ref{l33}, we now arrive at
\begin{eqnarray*}
  \|(R_{\delta} z_1 - R_{\delta} z_2)(t)\|_X &\leq& C \, \int_0^t \|z_1(\tau) -  z_2(\tau)\|_X d\tau,
\end{eqnarray*}
where $C$ depends on $\beta, \alpha, T, \|LG\|, \|G^*L^*\|$ and $\| (\delta I + LG G^*L^*)^{-1}\|$ and hence,
$$ \|R_{\delta} z_1 - R_{\delta} z_2 \|_Z \leq \frac{C T}{\sqrt{2}} \|z_1 - z_2\|_Z.$$

Proceeding inductively, we obtain that there exists a constant
$\gamma_n = \frac{(2 C T)^n}{\sqrt{2 n} (3 \cdot 5 \cdots 2 n -1)}$,
such that
$$ \|R^n_{\delta} z_1 - R^n_{\delta} z_2\|_Z \leq \gamma_n \|z_1 - z_2\|_Z.$$
Choose $n$ large enough (independent of $T$ and $C$)
such that $\gamma_n < 1,$ and hence, $R_{\delta}^n$ is a contraction. Therefore, by Banach contraction
mapping theorem, $R_{\delta}^n$ has a unique fixed point,say, $y_{\delta}^*$,
which is the limit of the sequence defined by (\ref{e3.9}). This $y_{\delta}^*$ is also the unique fixed
point of the operator $R_{\delta},$ for fixed $\delta \in (0,1].$

Next to show that $\mathcal{M} z_{\delta,k} \rightarrow \mathcal{M} y_{\delta}^*$. Setting $u_{\delta,k} =
\mathcal{M} z_{\delta,k},$
where $z_{\delta,k}$ is the mild solution of the system (\ref{e3.2}) with control $u_{\delta,k}.$ Then, we obtain
\begin{eqnarray*}
\left\|\Big(\mathcal{M} z_{\delta,k} - \mathcal{M} y_{\delta}^*\Big)(t)\right\|_X &\leq& C T^{1/2}
\|G^*L^* (\epsilon I + LG G^*L^*)^{-1}\| \;
\|z_{\delta,k} - y^*_{\delta}\|_Z,
\end{eqnarray*}
and hence,
\begin{eqnarray*}
\left\|\Big(\mathcal{M} z_{\delta,k} - \mathcal{M} y^*_{\delta}\Big) \right\|_Z &\leq& C T
\|G^*L^* (\epsilon I + LG G^*L^*)^{-1}\|\; \|z_{\delta,k} -y^*_{\delta}\|_Z.
\end{eqnarray*}
Since for each fixed $\delta \in (0,1],$ the sequence $z_{\delta,k} \rightarrow y^*_{\delta}$ in $Z$, this implies that $\mathcal{M} z_{\delta,k}
\rightarrow \mathcal{M} y^*_{\delta} = u^*_{\delta}.$

Since $y^*_{\delta}$ is the mild solution of the system (\ref{e1.1}) with control $u_{\delta}^*$. As $z_{\delta,k}
\rightarrow y^*_{\delta}$, it
follows that $R_{\delta} z_{\delta,k} \rightarrow R_{\delta} y^*_{\delta} = y^*_{\delta}.$
Using  the definition of $R_{\delta}$ and with  similar arguments
as earlier, we find that
$$ R_{\delta} z_{\delta,k}(t) = S(t) y_0 + \int_0^t S(t - \tau) \tilde{B}z_{\delta,k}(\tau)\; d\tau
+ \int_0^t S(t - \tau) Gu_{\delta,k}(\tau)\; d\tau.$$
As $k \rightarrow \infty,$ we obtain
$$ y^*_{\delta}(t) = S(t) y_0 + \int_0^t S(t - \tau) \tilde{B} y^*_{\delta}(\tau)\; d\tau
+ \int_0^t S(t - \tau) Gu^*_{\delta}(\tau)\; d\tau, $$
and $y^*_{\delta}$ is the mild solution of the system (\ref{e1.1}),
corresponding to control $u^*_{\delta}$ given by
\begin{eqnarray}\label{eq:u*-delta}
u^*_{\delta} &=&\mathcal{M} y_{\delta}^*\\
&=& G^*L^*\Big(\delta I + LG G^*L^*\Big)^{-1} \left[\hat{z}- L\tilde{B} y^*_{\delta}\right].
\end{eqnarray}

It remains to show that the problem (\ref{e1.1}) is approximately controllable.
To this end observe that
\begin{eqnarray}\label{eq:lim}
LG u^*_{\delta} &=& LG G^*L^*\Big(\delta I + LG G^*L^*\Big)^{-1} \left[\hat{z}- L\tilde{B} y^*_{\delta}\right] \nonumber\\
&=& \Big( (\delta I +LG G^*L^*)- \delta I\Big)\;\Big(\delta I
+ LG G^*L^*\Big)^{-1} \left[\hat{z}- L\tilde{B} y^*_{\delta}\right] \nonumber\\
&=& [\hat{z}- L\tilde{B} y^*_{\delta}] - \delta\;\Big(\delta I
+ LG G^*L^*\Big)^{-1} \left[\hat{z}- L\tilde{B} y^*_{\delta}\right].
\end{eqnarray}
Since $\|\hat{z}- L\tilde{B} y^*_{\delta}\|$ is bounded, a use of  Lemma \ref{l3.1} (iv) yields
$$\lim_{\delta \to 0^{+}}\left\|- \delta\;\Big(\delta I
+ LG G^*L^*\Big)^{-1} \left[\hat{z}- L\tilde{B} y^*_{\delta}\right]\right\| =0,$$
and hence,
$$\lim_{\delta \to 0^{+}} \| LG u^*_{\delta}+ L\tilde{B} y^*_{\delta}-\hat{z}\| =0,$$
that is, for given any given $\delta_1 >0,$ there exists a $\delta_0 >0$ such that for $ 0< \delta \leq \delta_0$
 $$  \| LG u^*_{\delta}+ L\tilde{B} y^*_{\delta}-\hat{z}\| < \delta_1.$$
 Hence, the system (\ref{e1.1}) is approximately controllable.
This completes the rest of the proof.
\hfill{$\blacksquare$}



\begin{remark}
Note that $U_{\delta} \neq \emptyset.$
Further, the  error of  the approximation in this case is given by
$$e_{\delta} z^{*}_{\delta} = \delta \;\Big(\delta I
+ LG G^*L^*\Big)^{-1}\left[\hat{z}- L\tilde{B} y^*_{\delta}\right].$$
\end{remark}
\begin{remark}
Under assumption $(A1)$-$(A4),$ the Theorem \ref{t11} implies that the system (\ref{e1.1}) is controllable without any
inequality constraint on $T$.
\end{remark}

Thus, we have $U_{\delta} \neq \emptyset$. The pair $(u^*_{\delta}, y^*_{\delta})$ so
obtained need not be an
optimal pair satisfying (\ref{e1.6}), and hence, the  problem (\ref{e1.6}) remains unanswered.


We now change our strategy and examine the process of obtaining the optimal pair of the constrained problem
through a sequence of optimal pairs of the unconstrained problems, as indicated in the Section 1. For this purpose, we first define a sequence of functionals $\{J_{\epsilon}\}$ with $\epsilon >0$ as
\begin{eqnarray}\label{e3.10}
J_{\epsilon}(u) = \frac{1}{2} J(u) + \frac{1}{2\epsilon} P(u), \; u \in Y,
\end{eqnarray}
where penalty function $P(u)$ is of the form
\begin{eqnarray}\label{e3.11}
P(u) = \left\|L Gu + L \tilde{B} W u - \hat{z} \right\|_X^2, \; u
\in Y.
\end{eqnarray}

Now the problem under investigation is to seek $u_{\epsilon}^* \in U$ such that
\begin{eqnarray}\label{e3.12}
J_{\epsilon} (u_{\epsilon}^*) = \inf_{u \in Y} J_{\epsilon}(u).
\end{eqnarray}

As in \cite{G-L-H}, roughly speaking, the approximate controllability can be viewed as the limit
of a sequence of optimal control problems (\ref{e3.12}).

We now make further assumption  that
\begin{description}
\item[{\bf $(A5)$}] The solution operator $W : U \rightarrow
Z$ is completely continuous.
\end{description}

\begin{remark}
One of the sufficient condition for $W$ to be completely continuous is that the semigroup
$\{S(t)\}$ is compact.
\end{remark}

Denote by $\mathcal{E}$ the operator $\mathcal{E} u = L Gu + L \tilde{B} W u$, where the operators
$L, \tilde{B}$ and $W$ are as defined before. Then, the functional $J_{\epsilon}$ defined
through (\ref{e3.10}) can be written as
\begin{eqnarray}\label{e3.13}
J_{\epsilon} (u) =  \frac{1}{2} \|u\|_Y^2 + \frac{1}{2\epsilon} \|\mathcal{E} u -\hat{z} \|_X^2.
\end{eqnarray}
Note that the operator $\mathcal{E}$ is a sum of linear continuous operator $L$ and a completely
continuous operator $W$ and hence, it is a weakly continuous operator.

\begin{theorem}\label{T3.2}
Under assumptions $(A1)$-$(A5),$ the unconstrained optimal control problem (\ref{e3.12}) has an optimal pair $(u_{\epsilon}^*,
y_{\epsilon}^*)$ such that $u_{\epsilon}^* \in U$ minimizes $J_{\epsilon}(u)$ and $y_{\epsilon}^*$
solves (\ref{e1.1}) corresponding to the control $u_{\epsilon}^*$.
\end{theorem}
{\bf Proof.} We first prove the weakly lower semicontinuity of the functional $J_{\epsilon}$. Let $u_{\epsilon}^{n} \rightharpoonup u_{\epsilon}^*$ in $Y$, then, it follows that
\begin{eqnarray*}
\liminf_{n \rightarrow \infty} J_{\epsilon}(u_{\epsilon}^{n}) &=&  \liminf_{n
\rightarrow \infty} \left[\frac{1}{2} \| u_{\epsilon}^{n}\|_Y^2 +
\frac{1}{2\epsilon} \|\mathcal{E} u_{\epsilon}^{n} - \hat{z} \|_X^2 \right] \\
&\geq& \liminf_{n \rightarrow \infty} \frac{1}{2} \| u_{\epsilon}^{n}
\|_Y^2 + \liminf_{n \rightarrow \infty} \frac{1}{2\epsilon} \| \mathcal{E}
u_{\epsilon}^{n} - \hat{z} \|_X^2.
\end{eqnarray*}
Observe that
\begin{eqnarray*}
\| \mathcal{E} u_{\epsilon}^{n} - \hat{z} \|_X^2 = \|L Gu_{\epsilon}^{n}\|_X^2 + \|L \tilde{B} W
u_{\epsilon}^{n} - \hat{z}\|_X^2 + 2 \left< L u_{\epsilon}^{n}, L \tilde{B} W
u_{\epsilon}^{n}-\hat{z} \right>_X,
\end{eqnarray*}
and hence,
\begin{eqnarray*}
\liminf_{n \rightarrow \infty}  \| \mathcal{E} u_{\epsilon}^{n} - \hat{z} \|_X^2
&\geq& \liminf_{n \rightarrow \infty} \|LG u_{\epsilon}^{n}\|_X^2 +
\liminf_{n \rightarrow \infty} \|L \tilde{B} W u_{\epsilon}^{n} - \hat{z}\|_X^2 \\
&& ~~~~~~+ 2 \liminf_{n \rightarrow \infty} \left< LG u_{\epsilon}^{n}, L
\tilde{B} W u_{\epsilon}^{n} - \hat{z} \right>_X.
\end{eqnarray*}
From Lemma \ref{l2}, we  arrive at
\begin{eqnarray*}
\|L \tilde{B} \left(W u_{\epsilon}^n - W u_{\epsilon}^*\right)\|_X &\leq& C \int_0^T
\|(W u_{\epsilon}^n - W u_{\epsilon}^*)\|_X ds \\
&\leq& C T^{1/2} \|W u_{\epsilon}^n - W u_{\epsilon}^*\|_Z.
\end{eqnarray*}
Since $u_{\epsilon}^{n} \rightharpoonup u_{\epsilon}^*$ in $Y$, and $W$ is completely
continuous, this implies that $W u_{\epsilon}^n \rightarrow W u_{\epsilon}^*$ in $Z$ and
hence, $L \tilde{B} W u_{\epsilon}^n - \hat{z} \rightarrow L \tilde{B} W u_{\epsilon}^* - \hat{z}$ in $X$. Using the fact that $L$ is weakly continuous and $W$ is completely
continuous, we obain $L Gu_{\epsilon}^{n} \rightharpoonup L Gu_{\epsilon}^*, \; L \tilde{B} W u_{\epsilon}^{n} - \hat{z} \rightarrow L \tilde{B} W u_{\epsilon}^* -
\hat{z}$ and $\left<L Gu_{\epsilon}^{n}, L \tilde{B} W {\epsilon}^{n} - \hat{z}\right> \rightarrow \left<LG u_{\epsilon}^*, L \tilde{B}W u_{\epsilon}^* -
\hat{z}\right>$ and along with the fact that the norm is weakly lower semicontinuous functional, we find that
\begin{eqnarray*}
\liminf_{n \rightarrow \infty} J_{\epsilon}(u_{\epsilon}^{n}) \geq \frac{1}{2}
\|u_{\epsilon}^*\|_Y^2 + \frac{1}{2\epsilon} \|\mathcal{E} u_{\epsilon}^* - \hat{z}\|_X^2.
\end{eqnarray*}
This proves the weakly lower semi-continuity of $J_{\epsilon}$.

Let $\{u_{\epsilon}^{n}\}$ be a minimizing sequence for the functional $J_{\epsilon}$, 
that is, $\inf_{u \in Y} J_{\epsilon}(u) = \lim_{n \rightarrow \infty}
J_{\epsilon}(u_{\epsilon}^n)$. Since $J_{\epsilon}$ is coercive, the sequence $\{u_{\epsilon}^{n}\}$ is bounded in $Y$. Then, there exists a subsequence which is also
denoted by $\{u_{\epsilon}^{n}\}$ such that $u_{\epsilon}^{n}\rightharpoonup u_{\epsilon}^*$ weakly in $Y$.  Since the functional (\ref{e3.13}) is weakly lower semicontinuous in $Y$, we  arrive at
$$ \inf_{u \in Y} J_{\epsilon}(u) = \lim_{n \rightarrow \infty} J_{\epsilon}(u_{\epsilon}^{n}) 
= \liminf_{n \rightarrow \infty} J_{\epsilon}(u_{\epsilon}^{n}) \geq J_{\epsilon}(u_{\epsilon}^*).$$
Therefore, we  obtain
$$ J_{\epsilon}(u_{\epsilon}^*) = \inf_{u \in Y} J_{\epsilon}(u).$$
As $y_{\epsilon}^n = W u_{\epsilon}^n$ and $u_{\epsilon}^n \rightharpoonup u_{\epsilon}^*$, the complete continuity of $W$ implies $y_{\epsilon}^n \rightarrow y_{\epsilon}^*$, where $y_{\epsilon}^* = W u_{\epsilon}^*$. Thus, $(u_{\epsilon}^*, y_{\epsilon}^*)$ is the optimal pair for the unconstrained optimal control problem (\ref{e3.12})
and this completes the proof of the theorem. \hfill{$\blacksquare$}

In our subsequent analysis, we need the following properties of the sequence
of minimizers $\{u_{\epsilon}^*\}$.
\begin{lemma}\label{L3.3}
Let $\epsilon >0$ be arbitrary and let $u_{\epsilon} \in Y$ be a minimizer of $J_{\epsilon}(u)$ in $Y$, where $J_{\epsilon}(u)$ as defined by (\ref{e3.10}). For $\epsilon' < \epsilon$, the followings holds:
\begin{description}
\item[(i)] $J_{\epsilon}(u_{\epsilon}) \leq J_{{\epsilon'}}(u_{{\epsilon'}}).$
\item[(ii)] $P(u_{\epsilon}) \geq P(u_{{\epsilon'}}).$
\item[(iii)] $J(u_{\epsilon}) \leq J(u_{{\epsilon'}}).$
\item[(iv)] $J(u_{\epsilon}) \leq J_{\epsilon}(u_{\epsilon}) \leq J(u^*) + \frac{\delta_0^2}{2 \epsilon}.$
\end{description}
\end{lemma}
{\bf Proof}. As $u_{\epsilon}$ minimizes $J_{\epsilon}$, it follows that
$$J_{\epsilon}(u_{\epsilon})=J(u_{\epsilon})+\frac{1}{2\epsilon}P(u_{\epsilon}) \leq J(u_{\epsilon'})+\frac{1}{2\epsilon}P(u_{\epsilon'}) \leq J(u_{\epsilon'})+\frac{1}{2\epsilon'}P(u_{\epsilon'})=J_{\epsilon'}(u_{\epsilon'}).$$
This proves $(i).$ For $(ii),$ let $u_{\epsilon}$ and $u_{\epsilon'}$ be the minimizers of $J_{\epsilon}$
and $J_{\epsilon'},$ respectively, then, we arrive at
$$J_{\epsilon}(u_{\epsilon})=J(u_{\epsilon})+\frac{1}{2\epsilon}P(u_{\epsilon}) \leq J(u_{\epsilon'})+\frac{1}{2\epsilon}P(u_{\epsilon'})$$
and
$$J_{\epsilon'}(u_{\epsilon'})=J(u_{\epsilon'})+\frac{1}{2\epsilon'}P(u_{\epsilon'}) \leq J(u_{\epsilon})+\frac{1}{2\epsilon'}P(u_{\epsilon})$$
On adding above inequalities, we get
$$P(u_{\epsilon}) \geq P(u_{{\epsilon'}}).$$
For $(iii),$ note that
$$J_{\epsilon}(u_{\epsilon})=J(u_{\epsilon})+\frac{1}{2\epsilon}P(u_{\epsilon}) \leq J(u_{\epsilon'})+\frac{1}{2\epsilon}P(u_{\epsilon'}).$$
Hence, using $(ii)$ it follows that
$$J(u_{\epsilon})-J(u_{\epsilon'})\leq \frac{1}{2\epsilon}\left(P(u_{\epsilon'})-P(u_{\epsilon})\right)
\leq 0 $$
and
$$J(u_{\epsilon}) \leq J(u_{\epsilon'}).$$
For $(iv),$ again we observe that
$$J(u_{\epsilon}) \leq J_{\epsilon}(u_{\epsilon})=J(u_{\epsilon}) 
+ \frac{1}{2\epsilon}P(u_{\epsilon}) \leq J(u^*) + \frac{1}{2\epsilon} P(u^*)\leq J(u^*)
 + \frac{\delta_0^2}{2 \epsilon}.$$
This completes the rest of the proof. \hfill{$\blacksquare$}\\


We are now in a position to state the main theorem of this article.

\begin{theorem}\label{T3.3}
Assume that for a fixed $\delta > 0, \;  U_{\delta} \neq \emptyset$  and assumptions $(A1)$-$(A4)$ hold. Let $(u_{\epsilon}^*, y_{\epsilon}^*)$ be an optimal pair of
the unconstrained problem (\ref{e3.12}). As $\epsilon \rightarrow 0$,  there exists a subsequence of $(u_{\epsilon}^*, y_{\epsilon}^*)$ converges to
$(u_{\delta}^*, y_{\delta}^*)$, where $(u_{\delta}^*, y_{\delta}^*)$ is an optimal pair of the constrained optimal control problem (\ref{e1.6}).
Furthermore, if  $U_{\delta}$ is a singleton then the entire sequence  $(u_{\epsilon}^*, y_{\epsilon}^*)$ converges to
$(u_{\delta}^*, y_{\delta}^*)$.
\end{theorem}
{\bf Proof.} As for a fixed $\delta > 0, \;  U_{\delta} \neq \emptyset$; the existence of the optimal pair $(u_{\epsilon}^*, y_{\epsilon}^*)$ to the
unconstrained problem (\ref{e3.12}) follows from the Theorem \ref{T3.2}. Let $u_{\epsilon'}^* \in U_{\delta}$. From Lemma  \ref{L3.3}, we have
$$J_{\epsilon}(u_{\epsilon}^*)  \leq J_{\epsilon'}(u_{\epsilon'}^*)  \; \mbox{for} \; \epsilon ' < \epsilon  $$

Thus, $\{J_{\epsilon}(u_{\epsilon}^*)\}$ is a  monotone decreasing sequence which is bounded below and hence it   converges. Similarly, $\{J(u_{\epsilon}^*)\}$ is also a convergent  sequence. Now $\frac{1}{\epsilon}P(u_{\epsilon}^*)$,
being the difference of two convergent sequence, also converges, which in turn, implies that $P(u_{\epsilon}^*)
\rightarrow 0$ as $\epsilon \rightarrow 0.$ Hence,
$$ \lim_{\epsilon\to 0}\|\mathcal{E} u_{\epsilon}^* - \hat{z} \|_X = 0.$$
Since $\{u_{\epsilon}^*\}$ is a uniformly bounded sequence in $Y$, it has a subsequence, again denoted by $\{u_{\epsilon}^*\}$ such that
 $u_{\epsilon}^* \rightharpoonup u_{\delta}^*$ in $Y.$
Weak continuity of $\mathcal{E}$ implies that $\mathcal{E} u_{\epsilon}^* \rightharpoonup  \mathcal{E} u_{\delta}^*.$
Hence $ \mathcal{E} u_{\delta}^* = \hat{z}$ and  $u_{\delta}^* \in   U_{\delta} $.
By the weak lower
semicontinuity of the norm functional and Lemma \ref{L3.3}, we arrive at
$$\|u_{\delta}^*\|_{Y} \leq \liminf_{{\epsilon} \rightarrow 0} \|u_{\epsilon}^* \|_{Y} \leq \limsup_{{\epsilon}
\rightarrow 0} \|u_{\epsilon}^* \|_{Y} \leq \|u_{\delta}^*\|_{Y},$$
and hence,
$$ \lim_{\epsilon \rightarrow 0} \|u_{\epsilon}^*\|_{Y} = \|u_{\delta}^*\|_{Y}.$$
This along with the weak convergence of $u_{\epsilon}^*$ to $u_{\delta}^*$, implies that
\[ u_{\epsilon}^* \rightarrow   u_{\delta}^*  \;  \mbox{as} \; \epsilon \rightarrow 0. \]

Again from Lemma \ref{L3.3} and weak lower semicontinuity of the norm functional, we obtain the inequality
$$ J(u_{\delta}^*) \leq \liminf_{{\epsilon} \rightarrow 0} J(u_{\epsilon}^*)
\leq J(\tilde{u}), \; \tilde{u} \in U_{\delta}.$$
This, in turn, implies that
$$ J(u_{\delta}^*) \leq J(\tilde{u}) \; \forall \;  \tilde{u} \in U_{\delta}.$$
Therefore, $(u_{\delta}^*, y_{\delta}^*)$ is the optimal pair for the constrained
optimal problem (\ref{e1.6}). It is also clear that if  $U_{\delta}$ is a singleton then the entire sequence
$  u_{\epsilon}^*$ converges to $u_{\delta}^* \; \mbox{in} \; Y.$

Next, we show the convergence of $y_{\epsilon}^*$ to $y_{\delta}^*$ in
$Z$. From (\ref{e1.3}) and Lemma \ref{l2}, we obtain
\begin{eqnarray*}
  \|y_{\epsilon}^*(t) - y_{\delta}^{*}(t)\|_X &\leq& C \int_0^t \|y_{\epsilon}^*(\tau) -
  y_{\delta}^{*}(\tau)\|_X d\tau + \beta \int_0^t \|u_{\epsilon}^*(\tau) -
  u_{\delta}^{*}(\tau)\|_X d\tau \\
&\leq& C \int_0^t \|y_{\epsilon}^*(\tau) - y_{\delta}^*(\tau)\|_X d\tau + \beta
T^{1/2} \|u_{\epsilon}^* - u_{\delta}^*\|_Y.
\end{eqnarray*}
Using Gronwall's lemma, we arrive that
\begin{eqnarray*}
\|y_{\epsilon}^*(t) - y_{\delta}^*(t)\|_X &\leq& \beta T^{1/2} \|u_{\epsilon}^* - u_{\delta}^*\|_Y \, e^{C T},
\end{eqnarray*}
and hence,
\begin{eqnarray}\label{ee3.14}
\|y_{\epsilon}^* - y_{\delta}^*\|_Z &\leq& \beta T e^{C T} \|u_{\epsilon}^* - u_{\delta}^*\|_Y.
\end{eqnarray}
Since $u_{\epsilon}^* \rightarrow u_{\delta}^*$ in $Y$, from (\ref{ee3.14}), we obtain $y_{\epsilon}^*
\rightarrow y_{\delta}^*$ in $Z$ as $\epsilon \rightarrow 0$. This completes
the rest of the proof. \hfill{$\blacksquare$}

\section{Approximation theorems}
\setcounter{equation}{0}
In our analysis, we are interested in the computation of the optimal control pair for the unconstrained problem.
We first begin by establishing some properties of the operator arising from the derivative of the
functional $J_{\epsilon}$, which is defined as follows:
\begin{eqnarray}\label{e4.1}
J_{\epsilon}(u) = \frac{1}{2} \|u\|_Y^2 + \frac{1}{2\epsilon} \|\mathcal{E} u -
\hat{z}\|_X^2,
\end{eqnarray}
where $\hat{z}$ is a fixed element in $X$. We first recall
the unconstrained optimal control problem
\begin{eqnarray}\label{e4.2}
J_{\epsilon}(u) = \inf_{v \in Y} J_{\epsilon}(v).
\end{eqnarray}

\begin{lemma}
The critical point of the functional $J_{\epsilon}$ is given by the
solution of the operator equation
\begin{eqnarray}\label{e44.1}
u + \frac{1}{\epsilon} \mathcal{K} \left(\mathcal{E} u - \hat{z}\right) = 0
\end{eqnarray}
where $\mathcal{K} = (LG + L \tilde{B} W)^*, \; \mathcal{E}u = (LG + L \tilde{B}W)
u$ and $y = W u$.
\end{lemma}
{\bf Proof:} We note that
\begin{eqnarray*}
J_{\epsilon}(u+h v) - J_{\epsilon}(u) &=& \frac{1}{2}\left<u +
h v, u + h v\right> \\
&& + \frac{1}{2\epsilon} \left<(LG + L\tilde{B}W)(u +
hv)-\hat{z}, (LG + L\tilde{B}W)(u + hv)-\hat{z}\right> \\
&& - \frac{1}{2} \left<u,u\right> - \frac{1}{2\epsilon}\left<(LG +
L\tilde{B}W)(u)-\hat{z}, (LG + L\tilde{B}W)(u)-\hat{z}\right> \\
&=& h \left<u, v\right> + \frac{h^2}{2}
\left<v, v\right> + \frac{h}{\epsilon} \left<LG u
+L\tilde{B}Wu - \hat{z}, (LG + L\tilde{B}W)(v)\right> \\
&& + \frac{h^2}{2\epsilon}
\left<(LG + L\tilde{B}W)(v), (LG + L\tilde{B}W)(v)\right>.
\end{eqnarray*}
Then, $J_{\epsilon}^{'}(u)$ is given by
\begin{eqnarray*}
J_{\epsilon}^{'}(u) v &=& \lim_{h \rightarrow 0} \frac{J_{\epsilon}(u+h v) -
J_{\epsilon}(u)}{h}\\
&=& \left<u, v\right> + \frac{1}{\epsilon} \left<LGu + L\tilde{B}Wu - \hat{z}, (LG +
L\tilde{B}W)(v)\right> \\
&=& \left<u, v\right> + \frac{1}{\epsilon} \left<(LG +
L\tilde{B}W)^*(LG + L\tilde{B}W)u - \hat{z}, v\right>,
\end{eqnarray*}
and hence,
$$ J_{\epsilon}^{'}(u) = u + \frac{1}{\epsilon} (LG +
L\tilde{B}W)^*(LG + L\tilde{B}W)u - \hat{z}).$$
If $u$ is a critical point of $J_{\epsilon}$, then it follows that
$$u + \frac{1}{\epsilon} \mathcal{K} (LGu + L\tilde{B}Wu - \hat{z})= 0,$$
where $\mathcal{K} = (LG + L \tilde{B} W)^*$. This concludes the proof.\hfill{$\blacksquare$}

Note that, in the literature, the operator equation (\ref{e44.1}) is known as the Hamerstein equation (see, Joshi {\em et. al.} \cite{M1})  Also note that the operator $\mathcal{K}$ is bounded linear operator. We first assume that the critical point of $J_{\epsilon}$ is the unique minimizer of $J_{\epsilon}$. Then the minimizing problem (\ref{e4.2}) is equivalent to the following solvability problem in the space $Y$:
\begin{eqnarray}\label{e4.11}
u + \frac{1}{\epsilon} \mathcal{K} E u = \hat{w},
\end{eqnarray}
where $\hat{w} = \frac{1}{\epsilon} \mathcal{K} \hat{z}$.
We now first begin approximating the main problem in the following way. Consider a family $\{X_m\}$  of finite dimensional subspaces
of $X$ such that
$$X_1 \subset X_2 \subset \ldots \subset X_m \ldots \subset X \; \mbox{with}
\; \overline{\bigcup_{m=1}^{\infty}} X_m = X.$$

Let $\{\phi_i\}_{i=1}^{\infty}$ be a basis for $X$. The approximating scheme for the space $Y = L^2(0, T; X)$ is then given by
the family of subspace $Y_m = L^2(0, T; X_m)$ such that
$$
Y_1 \subset Y_2 \subset \ldots \subset Y_m \ldots \subset Y \; \mbox{with}
\; \overline{\bigcup_{m=1}^{\infty}} Y_m = Y.
$$
Note that, the solution of the system (\ref{e1.1}) is given by $y(t) = \sum_{i=1}^{\infty} \alpha_i \phi_i$
with the control $u(t) = \sum_{i=1}^{\infty} \beta_i \phi_i$.\\

Let $P_m: X \rightarrow X_m$ be the projection given by
$$ P_m [y(t)] = \sum_{i=1}^{m} \alpha_i \phi_i, \; t \in [0, T],$$
where $X_m = \mbox{span} \{\phi_1, \phi_2, \ldots, \phi_m\}$. Then, this induces in a natural way the projection $\tilde{P}_m : Y \rightarrow Y_m$ given by
$$(\tilde{P}_m y)(t) = P_m y(t).$$

The projections $P_m$ and $\tilde{P}_m$ generate the approximating
operators $ \mathcal{K}_m$ and $\mathcal{E}_m$ defined by $\mathcal{K}_m
= \tilde{P}_m \mathcal{K}$ and $\mathcal{E}_m u  = P_m \mathcal{E} u$.
Then, the approximated minimization problem is stated as:
Find $u_m \in Y_m$ such that
\begin{eqnarray}\label{e4.5}
J_{\epsilon, m} (u_m) = \inf_{u \in Y_m} \left[ J_{\epsilon, m} (u) =
\frac{1}{2} \|\tilde{P}_m u\|_{Y_m}^2 + \frac{1}{2\epsilon} \|P_m \mathcal{E} P_m u -
P_m \hat{z}\|_{X_m}^2 \right].
\end{eqnarray}

As in the case of problem (\ref{e3.12}), one can show that the problem (\ref{e4.5}) has a solution $u_m \in Y_m$, and hence, its critical point satisfying the operator equation in the approximating space $Y_m$ as
\begin{eqnarray}\label{e4.6}
u_m + \frac{1}{\epsilon} \mathcal{K}_m \left(\mathcal{E}_m u_m - P_m \hat{z}\right) = 0.
\end{eqnarray}

Following theorem shows that the solution for the problem (\ref{e4.6}) is uniformly bounded in $Y_m$ and the approximating pair $(u_m^*, y_m^*)$ converges to $(u^*, y^*),$ where $(u^*, y^*)$ is an optimal pair of
the constrained problem (\ref{e1.6}).

\begin{theorem}\label{t4.1}
Let $U_{\delta} \neq \emptyset$ and $u_m^*$ be the solution to the problem (\ref{e4.5}). Then $\{u_m^*\}$ is uniformly bounded in $Y_m$. If in addition, $J_{\epsilon}$ possesses a unique minimizer in $Y$ which is also the only critical point of $J_{\epsilon}$, then (\ref{e4.5}) has an optimal pair $(u_m^*, y_m^*)$ which converges to $(u^*, y^*)$ in $Y,$ where $(u^*, y^*)$ is an optimal pair of the constrained problem (\ref{e1.6}).
\end{theorem}
{\bf Proof:} Existence of the optimal pair $(u_m^*, y_m^*)$ to the optimal control problem (\ref{e4.5}) follows from Theorem
\ref{T3.2}.

Let $u^* \in U_{\delta}$, then from the definition of $U_{\delta}$,  we have $\|\mathcal{E} u^* -
\hat{z}\|\leq \delta.$
Define $u_m^* = \tilde{P}_m u^*$. Then
\begin{eqnarray*}
\frac{1}{2} \| u_m^* \|^2 &\leq& J_{\epsilon, m} (u_m^*) = \frac{1}{2}
\|\tilde{P}_m u^*\|^2 + \frac{1}{2\epsilon} \|P_m \mathcal{E} u_m^* -
P_m \hat{z}\|^2\\
&\leq& \frac{1}{2} \|\tilde{P}_m u^*\|^2 + \frac{1}{2\epsilon} \|P_m\|^2
\|\mathcal{E} u_m^* - \hat{z}\|^2\\
&\leq & \frac{1}{2} \|\tilde{P}_m\|^2\; \|u^*\|^2 + \frac{1}{\epsilon} \|P_m\|^2
\Big(\|\mathcal{E} u_m^* - \mathcal{E} u^*\|^2 + \|\mathcal{E} u^*-\hat{z}\|^2\Big).
\end{eqnarray*}
Since $u_m^* = \tilde{P}_m u^* \rightarrow u^*$ and $\mathcal{E}$ is weakly continuous, we have $\mathcal{E} u_m^* \rightharpoonup \mathcal{E} u^*$. Hence, both the term on right
hand side is bounded. Therefore $\{u_m^*\}$ is uniformly bounded.

Since $\{u_m^*\}$ is uniformly bounded, it has a subsequence, still denoted by $u_m^*$, which converges weakly to $u^*$ in $Y$. Then from the weak lower semicontinuity of the norm functional and Lemma
\ref{L3.3}, we arrive at
$$ \|u^*\| \leq \liminf_{m \rightarrow \infty} \|u_m^*\| \leq \limsup_{m \rightarrow \infty} \|u_m^*\| \leq \|u^*\|. $$
This implies
$$ \lim_{m \rightarrow \infty} \|u_m^*\| = \|u^*\|.$$
Together with the fact that $u_m^* \rightharpoonup u^*$ in $Y$, we obtain
$$ u_m^* \rightarrow u^* \;\;\mbox{in} \; Y, \; \mbox{as} \; m
\rightarrow \infty.$$
As $y_m^* = W u_m^*$ and $u_m^* \rightarrow u^*$, then continuity of
the solution operator $W$ implies that $y_m^* \rightarrow y^* = W
u^*$. This now completes the rest of the theorem. \hfill{$\blacksquare$}

The next step is to discretize in the direction of $t$.
This leads to finite dimensional subspaces $Y_m^k$ of each fixed $Y_m$ as follows
$$ Y_m^k = \left\{ y_m^k \in {\cal P}_0 : y_m^k\left|_{[t_{l-1},t_l]}\right. =
y_m^l,\, t_0 = 0,\, t_k = 1,\, t_l = l \Delta t, \,\Delta t =
1/k,\, 1 \leq l \leq k \right\}$$ where ${\cal P}_0$ is the space
of piecewise constant polynomials. It is clear that $Y_m^k $
satisfies the following property
$$ Y_m^1 \subset Y_m^2 \subset \ldots Y_m^k \subset \ldots
\subset Y_m \;\;  \mbox{with} \;\; \overline{\bigcup_{k
=1}^{\infty} Y_m^k} = Y_m. $$

We denote by $Q_m^k$, the orthogonal
projection from $Y_m$ to $Y_m^k$. This induces the operators $\mathcal{K}_m^k = Q_m^k
\mathcal{K}_m$ and $\mathcal{E}_m^k u = Q_m^k \mathcal{E}_m u$.\\


For a fixed $m$, we approximate the minimization problem (\ref{e4.5}) by the following minimization problem in the finite dimensional subspace $Y_m^k$ of $Y_m$.\\

Find $u_m^k \in Y_m^k $ such that
\begin{eqnarray}\label{e4.7}
\Phi(u_m^k) = \inf_{u_m \in Y_{m}^k} \left[ J_{\epsilon, m}^k(u_m) =
\frac{1}{2} \|Q_m^k u_m\|_{Y_m^k}^2 + \frac{1}{2\epsilon} \; \| Q_m^k
\mathcal{E}_m u_m - \hat{z}_m \|_{X_m}^2 \right].
\end{eqnarray}

The unique minimizer of the problem (\ref{e4.7}) is given by the critical point of $\Phi$, which is equivalent to the following solvability problem in the space $Y_m^k $.
\begin{eqnarray}\label{e4.8}
u_m^k + \frac{1}{\epsilon} \mathcal{K}_m^k \left(\mathcal{E}_m^k u_m^k -
P_m \hat{z}\right) = 0.
\end{eqnarray}


On the lines of Theorem \ref{t4.1}, we have the following
theorem giving the convergence of the approximation optimal pair
$(u_m^k, y_m^k)$ as $k \rightarrow \infty$ with $m$ fixed.

\noindent
\begin{theorem}\label{t4.2}
Let $\{u_m^k\}$ be the solution of the problem (\ref{e4.8}). Then
the approximating optimal pair $(u_m^k,
y_m^k)$ converges to $(u_m^*, y_m^*)$ in $Y_m$.
\end{theorem}

\section{Application}
\setcounter{equation}{0}
Let $\Omega$ be a bounded domain in  $\mathbb{R}^d$ with smooth boundary $ \partial \Omega $. For fixed $ T > 0 $, let  $Q = (0,T)\times\Omega $ and $\Sigma =  (0,T)\times \partial \Omega$. Let $A$ be a second order uniformly elliptic differential operator given by (\ref{ee1.2}). Further, assume that the operator $B(t,s)$ is an unbounded partial differential operator of order $\beta \leq 2$ given by (\ref{ee1.3}).

\vspace{0.3cm}
\noindent Set $X = L^2(\Omega),\; V = H_0^1(\Omega),\; D(A) = H^2(\Omega) \cap H_0^1(\Omega)$ and $D(B)= H^2(\Omega)$. Then the weak formulation of the problem (\ref{e1.1}) is given by
\begin{eqnarray}\label{appl4-1}
\left(y_t, \phi\right) + A(y, \phi) &=& \int_0^t  B(t,s;y(s),\phi) ds + (u,\phi) \;\;  \forall \, \phi \in V, \; t \in [0, T]\\
y(0) &=& y_0, \nonumber
\end{eqnarray}
where $A(\cdot,\cdot) : H_0^1(\Omega) \times H_0^1(\Omega)
\rightarrow \mathbb{R}$ and $B(t,s;\cdot,\cdot) : H_0^1(\Omega)
\times H_0^1(\Omega) \rightarrow \mathbb{R}$ are the continuous
bilinear forms corresponding the operators $A$ and $B(t,s)$
respectively, that is
\begin{eqnarray*}
A(y, \phi) &=& \int_{\Omega} \left(\sum_{i,j=1}^d a_{i j}(x)
\frac{\partial y}{\partial x_j} \frac{\partial \phi}{\partial x_i}
+ c(x) y \phi \right)dx,
\end{eqnarray*}
and
\begin{eqnarray*}
 B(t,s;y,\phi) &=& \int_{\Omega}
\left(\sum_{i,j=1}^d b_{ij}(t,s;x) \frac{\partial y}{\partial x_j}
\frac{\partial \phi}{\partial x_i} + \sum_{j=1}^d
b_j(t,s;x)\frac{\partial y}{\partial x_j}\phi + b_0(t,s;x) y \phi
\right)dx.
\end{eqnarray*}

Here, $(a_{ij}(x))$ is a symmetric and positive definite
matrix for all $x \in \bar{\Omega}$ and $a_0(x) \geq 0$. From Lumer--Phillips theorem (see, Pazy \cite{P}), $(-A)$ generates a $C_0$--semigroup. For $y_0 \in D(A)$,  the unique mild solution for the system (\ref{e1.1}) is given by
\begin{eqnarray}\label{appl4-2}
y(t) = S(t) y_0 + \int_0^t S(t - \tau) \tilde{B} y(s) ds +
\int_0^t S(t - \tau) u(s) ds.
\end{eqnarray}

\noindent For final time $t = T$, we obtain
\begin{eqnarray}\label{appl4-3}
y(T) = S(T) y_0 + L \tilde{B} y  + L u
\end{eqnarray}
where the operator $\tilde{B}$ and $L$ are defined as before.\\

Since all the hypotheses (A1-A4) are satisfied, an appeal to Theorem \ref{t11} ensures the approximate controllability
of (\ref{appl4-1}). Also, set $U = L^2(\Omega)$ and $Y=L^2(0,T,U)$, the solution operator $W: U \rightarrow Y$ is compact and an application to Theorem \ref{T3.2} and \ref{T3.3} shows the existence of optimal control.\\

Let $\{{\cal J}_h\}$ be a family of regular
triangulation of $\Omega$ with $0 < h < 1$. For $K \in {\cal
J}_h$, set $h_K = \mbox{diam}(K)$ and $h = \max(h_K)$. Let
$$V_h = \left\{ v_h \in C^0(\bar{\Omega}): v_h|_{K} \in {\cal P}_1(K),
K \in {\cal J}_h, v_h = 0 \; \mbox{on} \; \partial \Omega
\right\},$$ where ${\cal P}_1(K)$ is the space of linear
polynomials on $K$. Then, the semidiscrete Galerkin approximation
of (\ref{appl4-1}) is defined by
\begin{eqnarray}\label{appl4-4}
\left(y_{h,t}, \chi\right) + A(y_h, \chi) &=& \int_0^t
B(t,s;y_h(s),\chi) ds + (u,\chi) \;\;  \forall \, \chi \in V_h, \; t \in [0, T] \\
y_h(0) &=& y_{0 h}, \nonumber
\end{eqnarray}
where, $y_{0 h}$ is the approximation of $y_0$ in $V_h$.\\

Let $\{\varphi_i\}_{i=1}^{N_h}$ be a bases of the
finite element space $V_h$. Since $y_h(t) \in V_h$, we write
$$y_h(t) = \sum_{i =1}^{N_h} \alpha_i(t) \varphi_i(x),$$
where $\{\alpha_i\}_{i=1}^{N_h}$ satisfies
\begin{eqnarray}\label{appl4-44}
\sum_{i =1}^{N_h}\left[\left(\varphi_i, \varphi_j\right)
\alpha_i^{'}(t) + A(\varphi_i, \varphi_j)\, \alpha_i(t) \right.
 &-&  \left. \int_0^t
B(\varphi_i, \varphi_j)\alpha_i(s) ds \right] \nonumber \\
& =&  \left(u(t), \varphi_j\right), \;\; j = 1,2,\cdots, N_h, \\
\alpha_i (0) &=& \gamma_i. \nonumber
\end{eqnarray}

Here, $\gamma_i$ is the coefficient of $\varphi_i(x)$ in
the representation of $y_{0h}$, that is, $y_{0h} = \sum_{i
=1}^{N_h} \gamma_i \varphi_i(x)$. This is the first order system of ordinary differential equations.\\


In matrix form, system (\ref{appl4-44}) can be written as follows
\begin{eqnarray}\label{appl4-444}
M \alpha^{'} + {\cal A} \alpha  - \int_0^t B \alpha(s) ds = U
\end{eqnarray}
where $M = [M_{ij}]$ with $M_{ij} = (\phi_i, \phi_j)$, ${\cal A} = [{\cal A}_{ij}]$ with ${\cal A}_{i j} = A(\phi_i, \phi_j), \; B = B_{i j}$ with $B_{i j} = B(\phi_i, \phi_j)$ and $U = [U_j]$ with $U_j = (u, \phi_j)$. Note that the system (\ref{appl4-444}) leads to a system of ordinary integro--differential equations and since
the mass matrix $M$ is invertible, the system (\ref{appl4-444}) is uniquely solvable in $C^1(0,T)$.\\

Let $P_h : V \rightarrow V_h $ be the
$L^2$--projection and let $\{S_h(t)\}$ denote the finite element analogue of $S(t)$, defined by the semidiscrete equation (\ref{appl4-4}) with $u = 0, \; B = 0$. This operator on $V_h$ may  be defined as the semigroup generated by the discrete analogue $A_h : V_h \rightarrow V_h$ of $A$, where
$$(A_h v, \chi) = A(v, \chi) \;\; \forall \; v, \chi \in
V_h.$$


Define the discrete analogue  $B_h =
B_h(t,s) : V_h \rightarrow V_h$ of $B = B(t,s)$ by
$$(B_h(t,s)v,\chi) = B(t,s; v,\chi) \;\; \forall \; v, \chi
\in V_h, \; 0 \leq s \leq t \leq T.$$

Now we write the semidiscrete problem (\ref{appl4-4}) in an abstract form
\begin{eqnarray}\label{appl4-6}
y_{h,t} + A_h y_h &=& \int_0^t B_h(t,s)y_h(s) ds + P_h u \equiv
\tilde{B}_h y_h + P_h u, \;\; \mbox{for}\; t \in [0,T],\\
y_h(0) &=& P_h y_0. \nonumber
\end{eqnarray}

Using Duhamel's principle, the solution $y_h$ of the semidiscrete problem (\ref{appl4-6}) may be written as
\begin{eqnarray}\label{appl4-7}
y_h(t) = S_h(t) P_h y_{0} + \int_0^t S_h(t - s) \tilde{B}_h y_h(s)
ds + \int_0^t S_h(t - s) P_h u(s) ds.
\end{eqnarray}
At time $t=T$, equation (\ref{appl4-7}) becomes
\begin{eqnarray}\label{appl4-8}
y_h(T) = S_h(T) P_h y_{0} + L_h \tilde{B}_h y_h + L_h P_h u,
\end{eqnarray}
where $L_h$ is defined by
\begin{eqnarray}\label{appl4-9}
L_h v_h = \int_0^T S_h(T - \tau) v_h(\tau) d\tau.
\end{eqnarray}

\noindent
Setting $e = y_h-y$, we have
\begin{eqnarray}\label{app-10}
e(T) &=& \left(S_h(T)P_h-S(T)\right)y_0 \nonumber\\
&&~~~~~~~~~~~+\left(\int_0^TS_h(T-s)\tilde{B}_hy_h(s) ds-\int_0^T S(T-s) \tilde{B}y(s) ds\right) \nonumber \\
&& ~~~~~~~~~~~~~~+ \left(\int_0^T S_h(T-s)P_hu(s) ds - \int_0^T S(T-s)u(s) \nonumber ds\right)\\
&=& F_h(T) y_0 +\int_0^T F_h(T-s)\tilde{B}y(s) ds \nonumber\\
&& ~~~~~~ + \int_0^T S_h(T-s) \left(\tilde{B}_hy_h(s)-P_h \tilde{B}y(s) \right)ds \nonumber + \int_0^T F_h(T-s)u(s) ds \nonumber\\
&=& I_1 + I_2+ I_3+I_4,
\end{eqnarray}
where the operator $F$ is defined as $F_h((t) = S_h(t)P_h -S(t)$. For $F_h$, it is well known that (see, Theorem 3.1 of Bromble {\em et. al.} \cite{b1})
\begin{eqnarray}\label{app-11}
\|F_h(t) v \| \leq c h^s t^{-(s-m)/2} |v|_m, \; 0 \leq m \leq s \leq 2.
\end{eqnarray}
For $I_1$,  using estimates (\ref{app-11}) for $v = y_0$, we get
\begin{eqnarray}\label{app-12}
I_1 \leq \|F_h(T) y_0\| \leq c h^s T^{-(s-m)/2} |y_0|_m, \; 0 \leq m \leq s \leq 2.
\end{eqnarray}
Now for $I_4$,
\begin{eqnarray*}
I_4 = \int_0^T F_h (T-s) u(s) ds &\leq& \int_0^T \|F_h(T-s)u(s) \| ds\\
&\leq& ch^s \int_0^T (T-s)^{-(s-m)/2} |u(s)|_m ds\\
&\leq& ch^s \left(\int_0^T (T-s)^{-(s-m)} ds \right)^{1/2} \left(\int_0^T |u(s)|_m^2 ds\right)^{1/2}
\end{eqnarray*}

Here, take $0 < s-m < 1$ with $0 \leq m \leq s \leq 2$. For the estimates of $I_2$, a use of Lemma 4.3 in Zhang \cite{Z} (pp. 135) yields
\begin{eqnarray*}
\|I_2\| \leq \|\tilde{F}_h(\tilde{B})(T)\| &\leq& \|\int_0^T F_h(T-s) \tilde{B}y(s) ds \|\\
&\leq& ch^2 \int_0^T (T-s)^{-1/2}\|\tilde{B}y(s) \| ds\\
&\leq& ch^2 \int_0^T (T-s)^{-1/2} \left( \int_0^s \|B(T,\tau) y(\tau) \| d\tau \right) ds \\
&\leq& ch^2 \int_0^T (T-s)^{-1/2} \left( \int_0^s \|y(\tau)\|_2d\tau\right) ds\\
&\leq& ch^2 \int_0^T (T-s)^{-1/2} s^{1/2} \left(\int_0^T\|y(\tau)\|_2 d\tau \right) ds\\
&\leq& c(T) h^2 \int_0^T \|y(\tau)\|_2 d\tau\\
&\leq&  c(T) h^2 \left(\|y_0\| + \|u\|_{L^2(L^2)}\right)
\end{eqnarray*}

For the term $I_3$, we again follow the idea of the proof of \cite{Z},  that is, following the existence of $e_2$ term in (\cite{Z}, pp 135-138) to conclude the estimates of $I_3$ as
\begin{eqnarray}\label{app-13}
|S_h(t) \chi|_{q,h} \leq c t^{-(p-q)/2} |\chi |_{p,h}.
\end{eqnarray}

For $m  \leq 1$, we have
\begin{eqnarray*}
\left<I_3, \chi\right> &=& \left<\int_0^T \int_0^s S_h(T-s) (B_h y_h(\tau)-P_h B y(\tau))d\tau ds, \chi\right>\\
&=& \int_0^T \left<\int_0^s(B_h y_h(\tau)-P_h B y(\tau))d\tau, S_h(T-s) \chi\right> ds\\
&=& \int_0^T \int_0^s B(s, \tau, e(\tau), S_h(T-s) \chi) d\tau ds.
\end{eqnarray*}

Using (4.13) of \cite{Z} and (\ref{app-13}), we obtain
$$|\left<I_3, \chi\right>| \leq c(T) \int_0^T \|e(\tau)\|d\tau \|\chi\|,$$
and hence
$$\|I_3\|  = \sup_{0 \neq \chi \in L^2} \frac{|\left<I_3, \chi\right>|}{\|\chi\|} \leq c(T) \int_0^T \|e(T)\| d\tau.$$

On substitution in (\ref{app-10}), we arrive at
\begin{eqnarray*}
\|e(T)\| &\leq& c \{h^2T^{-1/2}\|y_0\| + h^{1-\delta_0}\left(\int_0^T (T-s)^{1-\delta_0}ds\right)^{1/2} \left(\int_0^T \|u(s)\|^2 ds\right)^{1/2}\\
&& ~~+ h^2\left(\|y_0\| + \|u\|_{L^2(L^2)} \right)\} + c(T) \int_0^T \|e(\tau)\|d\tau   \;\;\; (\mbox{for small} \; \delta_0).
\end{eqnarray*}

By Gronwall's lemma, we obtain for small $\delta_0$,
$$\|e(T)\| \leq c h^{1-\delta_0}\left(\|y_0\| +\|u\|_{L^2(L^2)}\right).$$

\noindent {\it {Full Discretization}}: Let $k$ be the step size in time, $t_n = n k,\; n = 0,1,2,\cdots, \; N = T/k$ and  let $u^n = u(t_n)$. For $\phi \in C[0, T]$ set
$$ \bar{\partial_t} \phi(t_n) = \frac{\phi(t_n) -
\phi(t_{n-1})}{\Delta t}$$

The approximation $y_h^n \in V_h$ of $y_h$ at time $t = t_n$ is now defined as a solution of
\begin{eqnarray*}\label{appl4-45}
\left(\bar{\partial_t}y_h^n, \chi\right) + A(y_h^n, \chi) &=&  Q^n(B(y_h, \chi)) + (u^n,\chi), \; \chi \in V_h, \; n = 1,2,\cdots,N, \\
y_h^0 &=& y_{0 h} \;\mbox{in} \; \Omega, \nonumber
\end{eqnarray*}
where we have used the left rectangular rule
$$Q^n(y) = \sum_{i=0}^{n-1} k y(t_j) \approx \int_o^{t_n} y(s) ds$$
to discretize the Volterra integral term.

As a consequence of Theorem~ \ref{t4.2}, it is possible to show that  an approximate pair $\{u_h^n, y_h^n\}$
converges to the optimal pair $\{u_h^*, y_h^*\}.$

\section{Numerical experiment}
\setcounter{equation}{0}
\noindent In this section, we present a numerical experiment to
illustrate the computation of the minimizer $u^*$. We consider the
following one dimensional initial--boundary value problem
\begin{eqnarray}\label{chapt4-num1}
\frac{\partial y}{\partial t} - \frac{\partial^2 y}{\partial x^2}
&=& \int_0^t B(t,s) y(s) ds + u(t,x), \;\;\, \mbox{on} \; (0, T) \times (0,1) \nonumber \\
y(t,x) &=& y_0(x) \;\;\;\; ~~~~~~~~~~~~~~~~~~~~~~~~~~~ x \in (0,1) \\
y(t,0) &=& 0 ~~=~~y(1,t) ~~~~~~~~~~~~~~~~~~~\, t \in [0, T]
\nonumber
\end{eqnarray}

\indent Set $T = 1, \; \Omega = (0, 1) \subset \mathbb{R}^1$  with $B(t,s) = \exp{(-\pi^2 (t -s))} I, \; y_0(x) = \sin(\pi x)$ and $\hat{y} = \exp(-\pi^2) \sin(\pi x)$. For this system $\hat{y} \in R(T,y_0)$, since $y(t,x) = \exp(-\pi^2 t) \sin(\pi x)$ is an exact solution of the system (\ref{chapt4-num1}) corresponding to the control function $u(t,x)  = -t \exp(-\pi^2 t)  \sin(\pi x)$ with $y(T,x) = \exp(-\pi^2) \sin(\pi x)$.

\vspace{0.3cm}
\indent Here, we choose $\Delta t, \; h$ and $N = 1/\Delta
t$. Using MOA algorithm (see, Joshi \emph{et. al.}  \cite{A2}), we compute $u^n, n =1,2,\ldots,N$ and
then plot the graph of numerical results for $N = 40$. In Figure
1, we plot the graph of the approximated state at time $T=1$ and
the given final state $\hat{y} = \exp(-\pi^2) \sin(\pi x)$
corresponding to the approximated optimal control $u^*$. The
approximated optimal control $u^*$ has been shown in Figure 2.
Figure 3(i) shows the surface of the computed state
corresponding to the optimal control $u^*$, whereas Figure 3(ii)
shows the surface of the exact solution of the system
(\ref{chapt4-num1}).\\

\begin{figure}[ht]
\begin{center}
\includegraphics*[width=9cm,height=7cm]{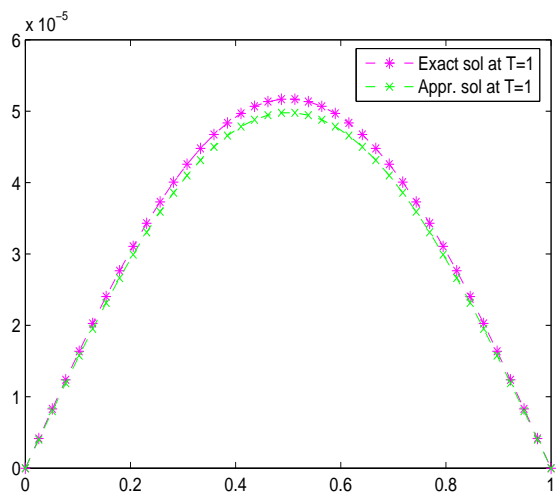}
\caption{Comparison between $y(T)$ and $\hat{y}$.}
\end{center}
\end{figure}

\begin{figure}[ht]
\begin{center}
\includegraphics*[width=9cm,height=7cm]{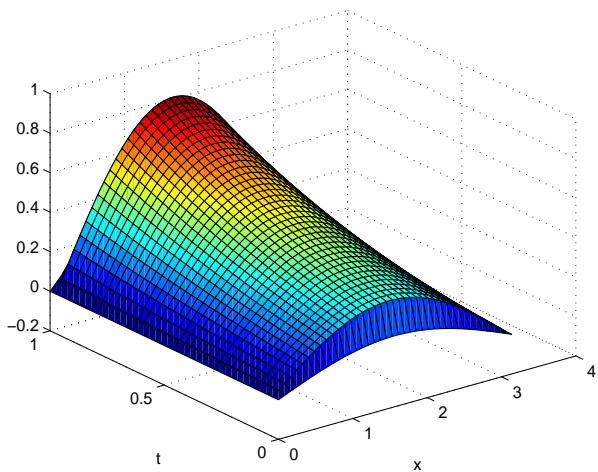}
\caption{Graph of the computed optimal control $u^*$.}
\end{center}
\end{figure}

\begin{figure}[ht]
\begin{center}
\includegraphics*[width=8cm,height=9cm]{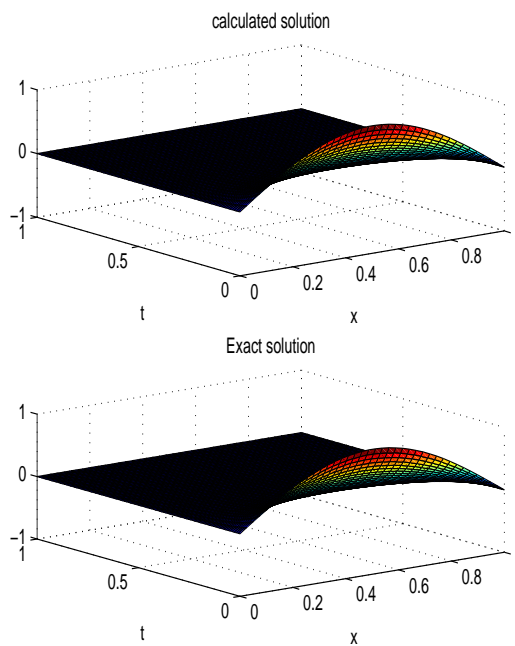}
\caption{Graph of (i) the computed state corresponding to the
optimal control $u^*$, and (ii) exact solution.}
\end{center}
\end{figure}

\bibliographystyle{plain}

\end{document}